\documentclass[reqno]{amsart}     
\usepackage{pb-diagram,enumerate,amsmath}     
\usepackage{pstricks,pst-node,pst-coil}     
\usepackage{graphicx}

\let\nc\newcommand 

\theoremstyle{plain}      
\newtheorem{step}{Step} 
\newtheorem{thm}{Theorem}[section]     
\newtheorem{theorem}[thm]{Theorem}

\newtheorem{lemma}[thm]{Lemma}     
\newtheorem{prop}[thm]{Proposition}     
\newtheorem{proposition}[thm]{Proposition}

\theoremstyle{remark}

\theoremstyle{definition}

\renewcommand{\epsilon}{\varepsilon}

\def\al{{\alpha}}

\def\De{{\Delta}}

\def\Si{{\Sigma}}         
         
\def\ep{{\varepsilon}}

\def\phi{{\varphi}}

\let\na\nabla     

\let\theta\vartheta
\let\phi\varphi

\DeclareMathAlphabet{\doba}{U}{msb}{m}{n}

\gdef\mN{\doba{N}}

\gdef\mR{\doba{R}}

\gdef\mZ{\doba{Z}}

\def\grad{{\mathop{\rm grad}}} 
\def\scal{{\mathop{\rm scal}}}

\def\eref#1{{\rm (\ref{#1})}}   

\let\<\langle 
\let\>\rangle

\nc{\PsoM}{\mbox{{$P_{\mbox{\scriptsize SO}}(M)$}}}
\nc{\PsonM}{\mbox{{$P_{\mbox{\scriptsize SO(n)}}(M)$}}}
\nc{\PsoG}{\mbox{{$P_{\mbox{\scriptsize SO}}(G)$}}}
\nc{\PspinM}{\mbox{{$P_{\mbox{\scriptsize Spin}}(M)$}}}
\nc{\PspinG}{\mbox{{$P_{\mbox{\scriptsize Spin}}(G)$}}}
\nc{\Pspineps}{\mbox{{$P_{\mbox{\scriptsize Spin,}\epsilon}$}}}

\begin{document}     

\title
[Surgery and Harmonic Spinors]
{Surgery and Harmonic Spinors}

\author{Bernd Ammann} 
\address{Institut \'Elie Cartan, BP 239 \\ 
Universit\'e de Nancy 1 \\
54506 Vandoeuvre-l\`es-Nancy Cedex \\ 
France}
\email{bernd.ammann@gmx.net}

\author{Mattias Dahl} 
\address{Institutionen f\"or Matematik \\
Kungliga Tekniska H\"ogskolan \\
100 44 Stockholm \\
Sweden}
\email{dahl@math.kth.se}

\author{Emmanuel Humbert} 
\address{Institut \'Elie Cartan, BP 239 \\ 
Universit\'e de Nancy 1 \\
54506 Vandoeuvre-l\`es-Nancy Cedex \\ 
France}
\email{humbert@iecn.u-nancy.fr}


\date{\today}


\thanks{This work was initiated when the authors visited the Albert
Einstein Institute, Golm, Germany, and was continued at the IECN Nancy, 
France. We wish to thank these institutes for their
very kind hospitality and support.}

\begin{abstract}
Let $M$ be a compact manifold with a fixed spin structure $\chi$.
The Atiyah-Singer index theorem implies that for any metric $g$ on $M$
the dimension of the kernel of the Dirac operator is bounded from below 
by a topological quantity depending only on $M$ and $\chi$.
We show that for generic metrics on $M$ this bound is attained.
\end{abstract}

\maketitle     

\section{Introduction}

We suppose that $M$ is a compact spin manifold.
By a \emph{spin manifold} we will always mean a smooth manifold equipped with 
an orientation and a spin structure.
After choosing a metric
$g$ on $M$, 
one can define the spinor bundle $\Si^gM$ and the Dirac operator
$D^g:\Gamma(\Si^g M)\to \Gamma(\Si^g M)$ 
see \cite{friedrich00,lawson.michelsohn89,hijazi01}. 

Being a self-adjoint elliptic operator
$D^g$ shares many properties with the Hodge-Laplacian 
$\De^g_p:\Gamma(\Lambda^p T^* M)\to \Gamma(\Lambda^p T^* M)$. 
In particular, if $M$ is compact, then the spectrum is discrete and real, and
the kernels of $\De^g_p$ and $D^g$ are finite-dimensional.
Elements of $\ker \De^g_p$ resp.\ $\ker D^g$ are called \emph{harmonic forms}
resp.\ \emph{harmonic spinors}.

However, the relation of $\Delta^g_p$ resp.\ $D^g$ to topology is different.
Hodge theory tells us that the Betti numbers $b_p:=\dim \ker \De^g_p$
only depend on the topological type of $M$. 
The dimension of the kernel of $D^g$ is invariant under conformal
changes of the metric, however it does depend on
the choice of conformal structure. The first examples
of this phenomenon  
were constructed by Hitchin \cite{hitchin74}, and 
it was conjectured by several people including B\"ar and 
the second named author \cite{baer_dahl02}
that $\dim \ker D^g$ depends on the metric for any compact spin manifold
of dimension $\geq 3$.  

On the other hand, $\dim \ker D^g$ is 
topologically obstructed. The Index
Theorem by Atiyah and Singer gives a topological lower bound on the
dimension of the kernel of the Dirac operator. For $M$ a compact spin
manifold of dimension $n$ this bound is \cite{lawson.michelsohn89},
\cite[Section~3]{baer_dahl02} 
\begin{equation} \label{boundfromAS}
\dim \ker D^g \geq 
\begin{cases}
|\hat{A}(M)|, &\text{if $n \equiv 0 \mod 4$;} \\
1, &\text{if $n \equiv 1 \mod 8$ and $\alpha(M)\neq 0$;} \\
2, &\text{if $n \equiv 2 \mod 8$ and $\alpha(M)\neq 0$;} \\
0, &\text{otherwise.} \\ 
\end{cases}
\end{equation}
Here the $\hat{A}$-genus $\hat{A}(M) \in \mZ$ and the $\alpha$-genus
$\alpha(M) \in \mZ/2\mZ$ are invariants of (the spin bordism class of)
the differential spin manifold $M$, and $g$ is any Riemannian metric on
$M$.

It is hence natural to ask whether metrics exist, such that 
equality holds in (\ref{boundfromAS}). Such metrics will be called 
\emph{$D$-minimal}.
In \cite{maier97} it is proved that a generic
metric on a manifold of dimension $\leq 4$ is $D$-minimal. In   
\cite{baer_dahl02} the same result is proved for manifolds of
dimension at least $5$ which are simply connected or have certain
fundamental groups. The argument in \cite{baer_dahl02} utilizes the
surgery-bordism method which has proven itself very powerful in the
study of manifolds with positive scalar curvature metrics. In a
similar fashion we will use surgery methods
to prove the following.
\begin{theorem} \label{main_generic}
Let $M$ be a compact connected spin manifold.
Then a generic metric on $M$ is $D$-minimal.
\end{theorem}
Our method also yields a new proof in dimensions $2$, $3$ and $4$.
Since $\dim \ker D$ behaves additively with respect to disjoint union
of spin manifolds while the $\hat{A}$-genus/$\alpha$-genus may cancel
it is easy to find disconnected manifolds with no $D$-minimal metric. 

Let us also mention that 
if $M$ is a compact Riemann surface of genus $\leq 2$, then \emph{all metrics} 
are $D$-minimal. The same holds for Riemann surface of genus $3$ whose 
spin structure is not spin bordant $0$. 
However if the genus is $\geq 4$ (or equal to $3$ with spin structures 
that are spin bordant $0$),
then there are also metrics with larger kernel \cite{hitchin74}, see also
\cite{baer.schmutz:92}.

In order to explain the surgery-bordism method in the proof of
Theorem~\ref{main_generic} we have to fix some notation.

A smooth embedding $f:N\to M$ is called \emph{spin preserving} if the pullback
of the orientation and spin structure of $M$ to $N$ under $f$ is 
the orientation and spin structure of $N$.
If $M$ is a spin manifold we
denote by $M^-$ the same manifold with the opposite orientation.

For $l \geq 1$ we denote by $B^l(R)$ the standard $l$-dimensional open ball
of radius $R$ and by $S^{l-1}(R)$ its boundary. We abbreviate 
$B^l=B^l(1)$ and $S^{l-1}=S^{l-1}(1)$. The standard Riemannian metrics
on $B^l(R)$ and $S^{l-1}(R)$ are denoted by $g^{\text{flat}}$ and
$g^{\text{round}}$. We equip $S^{l-1}(R)$ with the 
\emph{bounding spin structure}, i.e.\ the spin structure
obtained by restricting the unique spin structure on $B^l(R)$ (if
$l>2$ the spin structure on $S^{l-1}(R)$ is unique, if $l=2$ it is
not). 

Let $f:S^k\times \overline{B^{n-k}}\to M$ be a spin prerserving embedding, 
Then we define  
  $$\widetilde{M}= \left(M\setminus f(S^k\times B^{n-k})\right) \cup \left(\overline{B^{k+1}}\times S^{n-k-1}\right)/\sim$$
where $\sim$ identifies the boundary of $S^k\times S^{n-k-1}$ with 
$f(S^k\times S^{n-k-1})$. The toplogical space $\widetilde{M}$ carries a 
differential structure and a spin structure such that the inclusions 
$M\setminus f(S^k\times B^{n-k})\hookrightarrow \widetilde{M}$ and 
$ \overline{B^{k+1}}\times S^{n-k-1}\hookrightarrow \widetilde{M}$
are spin preserving smooth embeddings. 

We say that $\widetilde{M}$ is obtained from $M$ \emph{by surgery of 
dimension $k$} or \emph{by surgery of codimension $n-k$}.

The proof of Theorem of Theorem~\ref{main_generic} relies on
the following surgery theorem.

\begin{theorem} \label{main}
Let $(M,g^M)$ be a compact $n$-dimensional Riemannian spin manifold.
Let $\widetilde{M}$ be obtained from $M$ by surgery in dimension $k$, 
$k \in \{0,1,\dots,n-2 \}$. Then $\widetilde{M}$ carries a metric 
$g^{\widetilde{M}}$ such that 
$$
\dim \ker D^{g^{\widetilde{M}}} \leq \dim \ker D^{g^M}.
$$
\end{theorem}

\section{Preliminaries}

\subsection{Spinor bundles for different metrics} \label{identification}

Let $M$ be a spin manifold of dimension $n$ and let $g, g'$ be
Riemannian metrics on $M$. The goal of this paragraph is to identify
the spinor bundles of $(M,g)$ and $(M,g')$ using the method of
Bourguignon and Gauduchon introduced in \cite{bourguignon_gauduchon92}. 

There exists a unique endomorphism $b^{g}_{g'}$ of $TM$ which is
positive, symmetric with respect to $g$, and satisfies 
$g(X,Y) = g'(b^{g}_{g'} X, b^{g}_{g'}Y)$ for all $X,Y \in TM$.
This endomorphism maps $g$-orthonormal frames at a point 
to $g'$-orthonormal frames at the same point and we get a map 
$b^{g}_{g'}: \mathrm{SO}(M,g) \to \mathrm{SO}(M,g')$
of $\mathrm{SO}(n)$-principal bundles. If we assume that
$\mathrm{Spin}(M,g)$ and $\mathrm{Spin}(M,g')$ are equivalent spin
structures on $M$ the map $b^{g}_{g'}$ lifts to a map $\beta^{g}_{g'}$
of $\mathrm{Spin}(n)$-principal bundles, 
\dgARROWLENGTH=2em   
$$
\begin{diagram}   
\node{\mathrm{Spin}(M,g)}\arrow[2]{e,t}
{\beta^{g}_{g'}}\arrow{s}\node[2]
{\mathrm{Spin}(M,g')}\arrow{s}\\   
\node{\mathrm{SO}(M,g)}\arrow[2]{e,t}
{b^{g}_{g'}}\node[2]{\mathrm{SO}(M,g')}   
\end{diagram}.
$$
From this we get a map between the spinor bundles $\Sigma^{g} M$ and
$\Sigma^{g'} M$ denoted by the same symbol and defined by   
\begin{eqnarray}\label{indentspinbun}   
\Sigma^{g} M = 
\mathrm{Spin}(M,g) \times_\sigma \Sigma_n 
&\to&
\mathrm{Spin}(M,g')\times_\sigma \Sigma_n\nonumber
= \Sigma^{g'} M \\   
\psi=[s,\phi)]
&\mapsto&
[\beta^{g}_{g'} s,\phi]  = \beta^{g}_{g'} \psi 
\end{eqnarray} 
where $(\sigma,\Sigma_n)$ is the complex spinor representation, and 
where $[s,\phi]$ denotes the equivalence class of 
$(s,\phi)\in \mathrm{Spin}(M,g) \times_\sigma \Sigma_n$ for the equivalence
relation given by the action of $\mathrm{Spin}(n)$. The map $\beta_{g'}^g$
preserves fiberwise length of spinors.
   
We define the Dirac operator $D^{g'}$ acting on sections of the spinor
bundle for $g$ by
$$
{}^{g\mkern-4mu}D^{g'}
= 
(\beta^{g}_{g'})^{-1} \circ D^{g'} \circ \beta^{g}_{g'}
$$
In \cite[Thm. 20]{bourguignon_gauduchon92} the operator 
${}^{g\mkern-4mu}D^{g'}$ is computed in terms of $D^g$ and some extra
terms which are small if $g$ and $g'$ are close. Formulated in a way
convenient for us the relationship is 
\begin{equation} \label{relD}
{}^{g\mkern-4mu}D^{g'} \psi 
= 
D^g \psi 
+
A^{g}_{g'}(\nabla^g \psi) 
+
B^{g}_{g'}(\psi) 
\end{equation}
where $A^{g}_{g'} \in \hom(T^*M \otimes \Sigma^{g} M, \Sigma^{g} M)$
satisfies 
\begin{equation} \label{boundA^g_g'}
| A^{g}_{g'} | \leq C | g - g' |_g
\end{equation}
and $B^{g}_{g'} \in \hom(\Sigma^{g} M, \Sigma^{g} M)$ satisfies 
\begin{equation} \label{boundB^g_g'}
| B^{g}_{g'} | \leq C ( | g - g' |_g + | \nabla^g (g - g') |_g ) 
\end{equation}
for some constant $C$.

In the special case that $g'$ and $g$ are conformal with $g'= F^2 g$
for a positive smooth function $F$ we have 
\begin{equation} \label{confD}
{}^{g\mkern-4mu}D^{g'}( F^{-\frac{n-1}{2}} \psi) 
= 
F^{-\frac{n+1}{2}} D^g \psi
\end{equation}
according to \cite{hitchin74,baum81,hijazi01}.

\subsection{Notations for spaces of spinors} \label{spaces}

Throughout the article $\phi$ and $\psi$ and its variants denote spinors,
i.e.\ sections of the spinor bundle. 
If $S$ is a closed or open subset 
of $M$, we write $C^k(S)$ both for the space of $k$ times
differentiable functions on $S$ and for the space of $k$ times differentiable 
spinors. As the bundle will be clear from the context, this will not 
lead to ambiguities. On $C^k(S)$ we define the norm
  $$\|\phi\|_{C^k(S)}:= \sum_{l=0}^k\sup_{x\in S} |\na^l\phi(x)|.$$
We sometimes write $\|\phi\|_{C^k(S,g)}$ instead of  $\|\phi\|_{C^k(S)}$ to 
indicate that the spinor bundle and the norm depend on $g$.
The analogous notation is used for Schauder spaces $C^{k,\al}$.

Similarly $L^2(S)=L^2(S,g)$ and $H_k^2(S)=H_k^2(S,g)$
denote the space of $L^2$-spinors and $H_k^2$-spinors. 
These spaces come with the norms
 $$\|\phi\|_{L^2(S,g)}^2:=\int_S |\phi|^2\,dv^g \qquad 
   \|\phi\|_{H_k^2(S,g)}^2:=\sum_{l=0}^k \int_S |\na^l\phi|^2\,dv^g.$$

Let $U$ be an open set. The set of locally $C^1$-spinors $C^1_{\text{loc}}(U)$ 
carries a topology such that $\phi_i\to\phi$ in $C^1_{\text{loc}}(U)$ 
if and only if $\phi_i\to\phi$ in $C^1(K)$ for any compact subset 
$K\subset U$.

\subsection{Regularity and elliptic estimates}

In the following section $M$ is not necessarily compact.

\begin{lemma} Let $(M,g)$ be a Riemannian manifold, and let 
$\psi$ be a spinor of regularity $L^2$. If $\psi$ is \emph{weakly harmonic},
i.e.\  
  $$\int_M\<\psi,D\phi\>\,dv^g=0$$
for all compactly supported smooth spinors $\phi$, 
then $\psi$ is smooth.
\end{lemma}

\begin{lemma}\label{lem.reg}
Let $(M,g)$ be a Riemannian manifold and let $K \subset M$ a 
compact subset. Then there is a constant $C = C(K,M,g)$ such that
$$
\| \psi \|_{C^2(K,g)} \leq C \| \psi \|_{L^2(M,g)}
$$
for all harmonic spinors $\psi$ on $(M,g)$. 
\end{lemma}

{\bf Proof of the lemmata.}\\
The condition of the first 
lemma implies $\int \<\psi,D^2\Phi\>\,dv^g =0$ for any
compactly supported smooth spinor $\Phi$. Writing down the equation in local
coordinates, one can use standard tools from partial differential equations
(as for example \cite[Theorem 8.13]{gilbarg.trudinger77}) to derive 
via recursion that $\psi$ is contained in $H_k^2(K_1)$ for any $k\in \mN$ and
any $K_1$ compact in $M$, and
that  
\begin{equation}\label{hk.est}
\|\psi\|_{H_k^2(K_1,g)}\leq C\|\psi\|_{L^2(M,g)}.
\end{equation}
Suppose that the boundary of $K_1$ is smooth.
One then uses the Sobolev embedding $H_k^2(K_1,g)\to C^1(K_1,g)$
for $k>n/2+1$ 
(see \cite[Theorem~6.2]{adams75}), and we get $\psi\in C^1(K_1,g)$ 
and an estimate for $\|\psi\|_{C^1(K_1,g)}$
analogous to~\eref{hk.est}. Now one can use 
Schauder estimates as in \cite[Theorem 6.6]{gilbarg.trudinger77}
to conclude that $\psi$ is smooth on any compactum $K$ contained in
the interior of $K_1$, and in order to derive a $C^2$ estimate.
\qed

\begin{lemma}[{Ascoli's theorem, \cite[Theorem 1.30 and 1.31]{adams75}}]\label{lem.ascoli}
Let $\phi_i$ be a sequence bounded in $C^{1,\alpha}(K)$. 
Then a subsequence converges in $C^1(K)$.
\end{lemma}

\subsection{Removal of singularities lemma}

In the proof of Theorem \ref{main} we will need the following lemma.

\begin{lemma} \label{lem.sing.rem}
Let $(M,g)$ be an $n$-dimensional 
Riemannian spin manifold and let $S \subset M$ be a 
compact submanifold of dimension $k \leq n-2$. 
Assume that $\phi$ is a spinor field 
such that $\|\phi\|_{L^2(M)} < \infty$ and $D^g \phi = 0$ weakly on 
$M \setminus S$. Then $D^g \phi = 0$ holds weakly also on $M$.
\end{lemma}

\begin{proof}
Let $\psi$ be a smooth spinor compactly supported in $M$. We have to 
show that 
\begin{equation} \label{eq.singsol.show}
\int_M \<\phi,D^g \psi\>  \,dv^g = 0.
\end{equation}

Let $U_{S}(\epsilon)$ be the set of points of distance at most 
$\epsilon$ to $S$. 
For a small $\epsilon>0$ we choose a smooth function 
$\eta : M \to [0,1]$ such that $\eta = 1$ on 
$U_S(\epsilon)$, $|\grad \eta| \leq 2/\epsilon$ and 
$\eta = 0$ outside $U_S(2\epsilon)$. We rewrite the left hand 
side of (\ref{eq.singsol.show}) as
\begin{equation*}
\begin{split}
\int_M \<\phi,D^g \psi\> \,dv^g
&=
\int_M \< \phi,D^g ( (1-\eta) \psi + \eta \psi ) \>
\,dv^g \\
&=
\int_M \< \phi,D^g ( (1-\eta) \psi ) \> \,dv^g \\
&
+ 
\int_M \< \phi,\eta D^g \psi \> \,dv^g
+ 
\int_M \< \phi,\grad \eta \cdot \psi \> \,dv^g. 
\end{split}
\end{equation*}
As $D^g \phi = 0$ weakly on $M \setminus S$ the first term vanishes. 
The absolute value of the second 
term is bounded by 
$$
\| \phi \|_{L^2(U_S(2\epsilon))} \| D^g\psi \|_{L^2(U_S(2\epsilon))}
$$
which tends to $0$ as $\epsilon \to 0$. Finally, the absolute value of 
the third term is bounded by
\begin{equation*}
\begin{split}
\frac{2}{\epsilon} \| \phi \|_{L^2(U_S(2\epsilon))}  
\| \psi \|_{L^2(U_S(2\epsilon))}
&\leq
\frac{C}{\epsilon} \|\phi\|_{L^2(U_S(2\epsilon))}
{(\operatorname{Vol}(
U_S(2\epsilon) \cap \operatorname{supp}(\psi) 
)}^{\frac{1}{2}} \\
&\leq
C \|\phi\|_{L^2(U_S(2\epsilon))} \epsilon^{ \frac{n-k}{2} - 1}.
\end{split}
\end{equation*}
Since $n-k \geq 2$,
the third term also tends to $0$ as $\epsilon \to 0$.
\end{proof}

\subsection{Products with spheres}

The spectrum of $(D^{g^{\text{round}}})^2$ is bounded from below by
$l^2/4$.  

If $(M,g)$ and $(N,h)$ are compact Riemannian spin manifolds 
then the squared Dirac operator $(D^{g+h})^2$ on 
$(M \times N, g + h)$ can be identified with $(D^g)^2 + (D^h)^2$. 
We conclude the following.  

\begin{proposition} \label{appl.sbundles} 
Let $(M,g)$ be a compact spin manifold and $l \geq 1$. Then 
the spectrum of $(D^{g + g^{\text{round}}})^2$ on $M \times S^l$ is 
bounded from below by $l^2/4$.
\end{proposition}

\section{Proof of Theorem \ref{main}}

Our standing assumptions are: $(M,g)$ is a compact Riemannian spin 
manifold of dimension $n$ together with a $k$-dimensional submanifold~$S$ 
of~$M$ diffeomorphic to~$S^k$. We assume $n-k\geq 2$.
The restriction of $g$ to~$S$ is 
denoted by~$h$. Let $\nu \to S$ be the normal bundle of $S$. We assume 
furthermore that a trivialization of the normal bundle is given, that 
is a vector bundle map $\iota: \mR^{n-k}\times S \to \nu$. 
We assume that $\iota$ is fiberwise an isometry.

For $R > 0$ we denote by $\nu(R)$ the disk bundle of vectors of length 
$\leq R$ in $\nu$. For sufficiently small $R$ the normal exponential map 
$\exp^\nu$ of $S$ defines a diffeomorphism of $\nu(R)$ onto a neighborhood 
of $S$. For such small $R>0$ one has
$$
U_S(R)= (\exp^\nu \circ \iota) (\overline{B^{n-k}(R)}\times S)=\exp^\nu(\nu(R)). 
$$

\begin{lemma}\label{lem.ch.spin}
Let $n\geq 3$.
Let $\chi$ be the canonical spin structure on $\mR^{n-1}$,
let $\chi_b$ be the bounding spin structure on~$S^1$ and $\chi_{nb}$
the non-bounding spin structure on $S^1$. 
There is a diffeomorphism from $F:\mR^{n-1}\times S^1$ to itself preserving
the linear structure of $\mR^{n-1}$ with
  $$F^*(\chi\times \chi_b)= \chi\times \chi_{nb}.$$
\end{lemma}

\begin{proof}
Let $\gamma:S^1\to \mathrm{SO}(n-1)$ be a generator 
of~$\pi_1(\mathrm{SO}(n-1))$. Then the map $(X,x)\mapsto (\gamma(x)X,x)$
is a diffeomorphism as desired.
\end{proof}

Let $\exp^\nu:\nu\to M$ be the restriction of the exponential map to $\nu$.
Close to the zero section of $\nu$, $exp^\nu$ is a diffeomorphism onto 
its image, and hence for small $\ep>0$ the map
  $$I_\iota:\mR^{n-k}\times S,\qquad 
     (X,x)\mapsto \exp\left(R\frac{\iota(X,x)}{\sqrt{1+\|X\|^2}}\right)$$
is a diffeomorpism onto the interior of $U_S(R)$.
The spin structure on $M$ induces a spin structure on $\mR^{n-k}\times S$.
If $k\geq 2$, then the spin structure on $\mR^{n-k}\times S$ is
unique. However, in the case $k=1$, the induced spin structure
might be $\chi\times \chi_b$ or $\chi\times \chi_{nb}$. 
If the induced spin structure is $\chi\times \chi_{nb}$, we replace
$\iota$ by $\iota'=\iota\circ F$, and the spin structure induced by
$I_{\iota'}$ is $\chi \times \chi_b$. Hence, we can assume from now on
without loss of generality that the trivialization $\iota$ 
induces the spin structure $\chi\times \chi_b$.

\subsection{Approximation by a metric of product form near $S$}

In the following $r(x)$ denotes the distance from the point $x$ to $S$ 
with respect to the metric $g$. 

\begin{lemma} \label{LemmaTaylorMetric}
For sufficiently small $R>0$ there is a constant $C > 0$ so that
$$
G = g - ((\exp^\nu\circ \iota)^{-1})^* (g^{\text{flat}} + h)
$$
satisfies 
$$
|G(x)|\leq C r(x), \quad |\nabla G (x)| \leq C
$$
on $U_{S}(R)$.
\end{lemma}
Note that in this lemma the function $r(x)$ is by definition 
the distance of $x$ to $S$ with respect to $g$ but it coincides with 
the distance of $x$ to $S$ with respect to the metric 
$((\exp^\nu\circ \iota)^{-1})^* (g^{\text{flat}} + h)$
\begin{proof}
Since $x \mapsto \nabla G(x)$ is continuous on a neighborhood of $S$
we can find a constant $C$ such that $|\nabla G (x)| \leq C$ for 
sufficiently small $R>0$. Now, let $x\in S$.  At first
the spaces $T_xS$ and $\nu_x$ are orthogonal with respect to the 
two scalar products $g(x)$ and 
$((\exp^\nu\circ \iota)^{-1})^* (g^{\text{flat}} + h)(x)$. 
It is also clear that these two scalar products coincide on $T_xS$. 
Since the differential $d(\exp^\nu \circ \iota)$ is an isometry,
they coincide also on $\nu_x$. This implies that 
$g(x)= ((\exp^\nu\circ\iota)^{-1})^* (g^{\text{flat}} + h)(x)$ 
and hence that $G(x)=0$. We obtain that $G$ vanishes on $S$. 
Since $G$ is $C^1$, $|G|$ is $1$-lipschitzian and thus there 
exists $C>0$ such that $|G(x)| \leq C r(x)$. 
\end{proof}

The following proposition allows us to assume that the metric $g$ has 
product form close to the surgery sphere $S$.
\begin{proposition} \label{productform}
Let $(M,g)$ and $S$ be as above. Then there is a metric $\tilde g$ on
$M$ and $\epsilon>0$ such that $d^g(x,S) = d^{\tilde g}(x,S)$, 
$\tilde g$ has  product form on $U_S(\epsilon)$ and  
$$
\dim \ker D^{\tilde g} \leq \dim \ker D^g.
$$
\end{proposition}
For $\delta > 0$ let $\eta$ be a smooth cut-off function such that 
$0 \leq \eta \leq 1$, $\eta = 1$ on $U_S(\delta)$, $\eta = 0$ on 
$M \setminus U_S (2\delta)$, and $|d \eta|_g \leq 2/\delta$. We set  
$$
g_{\delta}
= 
\eta((\exp^\nu\circ \iota)^{-1})^* (g^{\text{flat}} + h) 
+ (1-\eta) g.
$$
Then $d^g(x,S) = d^{g_\delta}(x,S) = r(x)$. Through a series of lemmas 
we will prove the proposition for $\tilde g=g_{\delta}$ for $\delta$ 
sufficiently small. 

In the following estimates $C$ denotes a constant whose values
might vary from one line to another, which is independent of $\delta$ 
and $\eta$ but might depend on $M$, $g$, $S$. Terms denoted by 
$o_i(1)$ tend to zero when $i \to \infty$.
\begin{lemma}
Let $\delta_i$ be a sequence with $\delta_i \to 0$ as $i \to \infty$. 
Let $\phi_i$ be a sequence of spinors on $(M, g_{\delta_i})$ such 
that $D^{g_{\delta_i}} \phi_i = 0$ and 
$\int_M |\phi_i|^2 \,dv^{g_{\delta_i}} = 1$. Then the sequence 
$\beta^{g_{\delta_i}}_{g} \phi_i$ is bounded in $H_1^2(M,g)$.
\end{lemma}
\begin{proof}
As $\int |\beta^{g_{\delta_i}}_{g} \phi_i|^2 \,dv^g = 1 + o_i(1)$ we have
to show that 
$\alpha_i=\sqrt{\int_M |\nabla^g (\beta^{g_{\delta_i}}_{g}\phi_i)|_g^2
  \,dv^g}$ is bounded. We assume the opposite, that is 
$\alpha_i\to \infty$, and set 
$\psi_i= \alpha_i^{-1} \beta^{g_{\delta_i}}_{g} \phi_i$. 
Then we have 
${}^{g\mkern-4mu}D^{g_{\delta_i}} \psi_i = 0$
since $\beta^{g}_{g_{\delta_i}} \circ \beta^{g_{\delta_i}}_{g}
= \operatorname{Id}$, so formula (\ref{relD}) gives us  
\begin{eqnarray*}
1
&=& 
\int_M |\nabla^g \psi_i|_g^2 \,dv^g \\
&=&
\int_M (|D^g \psi_i|^2 
- \frac{1}{4} \scal^g |\psi_i|^2) \,dv^g \\
&=&
\int_M (
| A^{g}_{g_{\delta_i}}(\nabla^g \psi_i) 
+ B^{g}_{g_{\delta_i}}(\psi_i) |^2 
- \frac{1}{4} \scal^g |\psi_i|^2) \,dv^g \\
&\leq&
\int_M (
2 | A^{g}_{g_{\delta_i}}(\nabla^g \psi_i) |^2
+ 2 | B^{g}_{g_{\delta_i}}(\psi_i) |^2 
- \frac{1}{4} \scal^g |\psi_i|^2) \,dv^g. 
\end{eqnarray*}
Using (\ref{boundA^g_g'}), (\ref{boundB^g_g'}), Lemma
\ref{LemmaTaylorMetric}, and the fact that $g$ and $g_{\delta_i}$ 
coincide outside $U_S(2\delta_i)$ we get 
\begin{eqnarray*}           
1  
&\leq& 
C \delta_i^2 \int_{U_S(2\delta_i)} | \nabla^g \psi_i|_g^2 \,dv^g
+
C \int_{U_S(2\delta_i)} | \psi_i|^2 \,dv^g
+
C \int_M | \psi_i|^2 \,dv^g \\
&\leq&
C \delta_i^2
+
C \int_{U_S(2\delta_i)} | \psi_i|^2 \,dv^g
+
\alpha_i^{-2} (1 + o_i(1)) \\
&\leq&
C \int_{U_S(2\delta_i)} | \psi_i|^2 \,dv^g
+
o_i(1)
\end{eqnarray*}
As $\psi_i$ is bounded in $H_1^2(M,g)$, a subsequence converges weakly
in $H_1^2(M,g)$ and strongly in $L^2(M,g)$ to a limit spinor 
$\psi\in H_1^2(M,g)$. Hence for this subsequence
$$
\int_{U_{S}(2\delta_i)}|\psi_i|_g^2 \,dv^g \to 0
$$ 
which implies a contradiction.
\end{proof}
\begin{lemma} \label{LemmaConvHarmonic}
Again let $\delta_i$ be a sequence with $\delta_i \to 0$ as 
$i \to \infty$ and let $\phi_i$ be a sequence of spinors on 
$(M, g_{\delta_i})$ such that $D^{g_{\delta_i}} \phi_i = 0$ and 
$\int_M |\phi_i|^2 \,dv^{g_{\delta_i}} = 1$.
Then, after passing to a subsequence, 
$\beta^{g_{\delta_i}}_{g} \phi_i$ converges weakly in $H_1^2(M,g)$  
and strongly in $L^2(M,g)$ to a harmonic spinor on $(M,g)$.
\end{lemma}
\begin{proof}
According to the previous Lemma the sequence 
$\beta^{g_{\delta_i}}_{g} \phi_i$ is bounded in $H_1^2(M,g)$ and hence
a subsequence converges weakly in $H_1^2(M,g)$. After passing to a
subsequence once again we obtain strong convergence in
$L^2(M,g)$. Denote the limit spinor by $\phi$. 

For any $\epsilon>0$ Lemma \ref{lem.reg} implies that 
$\beta^{g_{\delta_i}}_{g} \phi_i$ is bounded in 
$C^2(M \setminus U_S(\epsilon) )$, and Lemma \ref{lem.ascoli} then
implies that a subsequence converges in 
$C^1(M \setminus U_S(\epsilon) )$. Hence the limit $\phi$ is in 
$C^1_{\text{loc}}(M\setminus S)$ and satisfies $D^g \phi = 0$ on 
$M\setminus U(S)$. Since $\phi$ is in $L^2(M,g)$ it follows from 
Lemma \ref{lem.sing.rem} that $\phi$ is a weak solution of 
$D\psi=0$ on $(M,g)$. By elliptic regularity theory $\phi$
is a strong solution and a harmonic spinor on $(M,g)$. 
\end{proof}
\begin{proof}[Proof of Proposition \ref{productform}]
Let $m=\liminf_{\delta\to 0} \dim \ker D^{g_{\delta}}$.
For sufficiently small $\delta$ let 
$\phi_{\delta}^1,\dots,\phi_{\delta}^m \in \ker D^{g_{\delta}}$ 
be spinors such that 
\begin{equation}\label{eq.orth}
\int_M \< \phi_{\delta}^j, \phi_{\delta}^k \> \,dv^{g_{\delta}}
=
\begin{cases}
1, &\text{if $j=k$;} \\
0, &\text{if $j\neq k$.}
\end{cases}
\end{equation}
According to Lemma \ref{LemmaConvHarmonic} there are spinors 
$\phi^1,\dots,\phi^m \in \ker D^g$ and a sequence $\delta_i\to 0$
such that $\beta^{g_{\delta_i}}_{g} \phi_{\delta_i}^j$
converges to $\phi^j$ weakly in $H_1^2(M,g)$ and 
strongly in $L^2(M,g)$ for $j=1,\dots,m$.
Because of strong $L^2$-convergence the orthogonality relation 
(\ref{eq.orth}) is preserved in the limit so $\dim \ker D^g\geq m$.
Hence there is a $\delta_0>0$ so that 
$\dim \ker D^{g_{\delta_0}} = m \leq \dim \ker D^g$ and the
Proposition is proved with $\tilde g = g_{\delta_0}$.  
\end{proof}

\subsection{Proof for metrics of product form near $S$}

We assume that $g$ is a product metric on $U_S(R_{\text{max}})$ 
for some $R_{\text{max}}>0$, as we may from Proposition \ref{productform}.   
In polar coordinates 
$(r,\Theta) \in (0,R_{\text{max}}) \times S^{n-k-1}$ on
$B^{n-k}(R_{\text{max}})$ we get
$$
g = g^{\text{flat}} + h = dr^2 + r^2 g^{\text{round}} +  h.
$$

\begin{figure}\label{fig.hierarchy}
\begin{center}
  $$0< \rho << r_0<r_1/2 << R_{\text{max}}$$
\end{center}
\caption{Hierachy of variables}
\end{figure}

Let $\rho >0$ be a small number which we will finally let tend to $0$ (see also Figure~\ref{fig.hierarchy}).
We decompose $M$ into three parts
\begin{enumerate}  
\item $M \setminus U_S(R_{\text{max}})$, 
\item $(\rho/2, R_{\text{max}}) \times S^{n-k-1}\times S^k$, 
\item $U_S(\rho/2)=B^{n-k}(\rho/2)\times S^k$.
\end{enumerate}
The manifold $\widetilde{M}$ is obtained by removing part (3) and by gluing in
$S^{n-k-1}\times B^{k+1}$, that is $\widetilde{M}$ is the union of 
\begin{enumerate}
\item $M \setminus U_S(R_{\text{max}})$, 
\item $(\rho/2, R_{\text{max}}) \times S^{n-k-1}\times S^k$, 
\item[(3')] $S^{n-k-1}\times B^{k+1}$.
\end{enumerate}
We now define a sequence of metrics $g_\rho$ on $\widetilde{M}$ 
such that the theorem holds for small $\rho>0$. The metrics $g_\rho$
will coincide with $g$ on part (1), but will be modified in part (2)
in order to close up nicely in part (3').
\begin{figure}\label{fig.picture}
\begin{center}
\includegraphics[scale=.8]{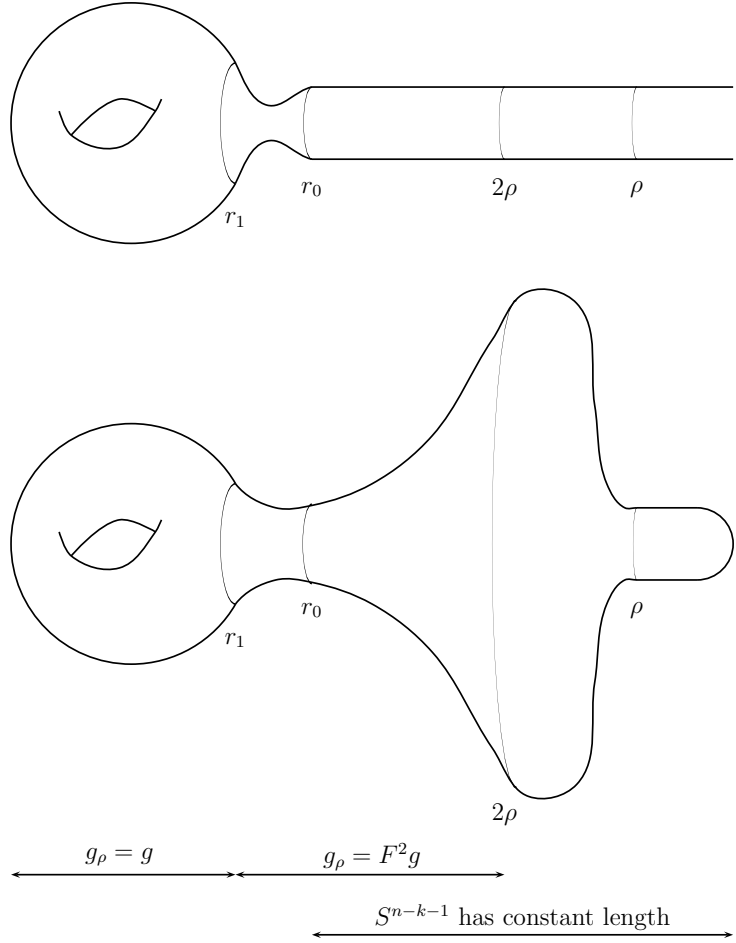}
\end{center}
\caption{The metric $g_\rho$.}
\end{figure}

Let $r_0, r_1$ be fixed such that 
$2 \rho < r_0 < r_1/2 < R_{\text{max}}/2$. 
Define $g_\rho$ on $\widetilde{M}$ by 
\begin{enumerate}
\item 
$g_{\rho} = g$ on $M \setminus U_S(R_{\text{max}})$, 
\item 
$g_{\rho} = F^2 ( dr^2 +r^2 g^{\text{round}}+ f^2_{\rho} h )$
on $(\rho/2, R_{\text{max}}) \times S^{n-k-1}\times S^k$, 
where $F$ and $f_{\rho}$ satisfy
$$
F(r)
=
\begin{cases}
1,   &\text{if $r_1 < r < R_{\text{max}} $;} \\
1/r, &\text{if $r < r_0$,}
\end{cases}
\quad
\text{and}
\quad
f_{\rho}(r)
=
\begin{cases}
1,  &\text{if $r > 2 \rho$;} \\
r,  &\text{if $r < \rho$ .}
\end{cases}
$$
\item [(3')]
$g_{\rho} = g^{\text{round}} + \gamma_{\rho}$ on
$S^{n-k-1} \times B^{k+1}$ where $\gamma_{\rho}$ is some 
metric so that $g_{\rho}$ is smooth. 
\end{enumerate}

The metric $g_\rho$ is visualized in Figure~\ref{fig.picture}.
In order to visualize the metric $g_\rho$ two projections are drawn.
In both projections the horizontal direction represents $-\log r$. 
In the first projection the vertical direction indicates the size
of the cosphere $S^{n-k-1}$. In the second projection the vertical 
direction indicates the size of $S$ which is fiberwise
homothetic to $(S\cong S^k,h)$.

We are now going to prove that
\begin{equation} \label{eqmain} 
\dim \ker D^{g_{\rho}} \leq \dim \ker D^g
\end{equation}
for small $\rho>0$. Before proving (\ref{eqmain}), we need some estimates.

For $\alpha \in (0,\rho/2)$, let $\widetilde{U}(\alpha) 
= \widetilde{M} \setminus (M \setminus U_S(\alpha))$ 
so that  $M \setminus U_S(\alpha) 
= \widetilde{M} \setminus \widetilde{U} (\alpha)$. 
\begin{prop} \label{estim1}
Let $s\in (0,r_1/2)$. Let $\psi_{\rho}$ be a harmonic spinor on
$(\widetilde{M},g_{\rho})$. Then for $\rho \in (0,s)$ it holds that
$$
\frac{(n-k-1)^2}{32} 
\int_{\widetilde{U}(s)\setminus \widetilde{U}(2 \rho)} 
|F^{\frac{n-1}{2}} \psi_{\rho}|^2 \,dv^g 
\leq  
\int_{\widetilde{U}(2s)\setminus \widetilde{U}(s)} 
|F^{\frac{n-1}{2}} \psi_{\rho}|^2 \,dv^g.  
$$
\end{prop}
\begin{proof}
Let $\eta \in C^{\infty}(\widetilde{M})$ be a cut-off function with 
$0 \leq \eta \leq 1$, $\eta = 1$ on $\widetilde{U}(s)$, $\eta = 0$ on 
$\widetilde{M} \setminus \widetilde{U}(2s)$, and 
\begin{equation} \label{deta}
|d\eta|_g \leq \frac{2}{s}.
\end{equation} 
The spinor $\eta \psi_\rho$ is compactly supported in $\widetilde{U}(2s)$. 
Moreover, the metric $g_\rho$ can be written as 
$g_\rho = g^{\text{round}} + h_{\rho}$ on $\widetilde{U}(2s)$ where 
the metric $h_\rho$ is equal to $r^{-2}dr^2 + r^{-2} f_{\rho}^2 h$ on 
$\widetilde{U}(2s)\setminus \widetilde{U}(\rho/2)$ and is equal to
$\gamma_{\rho}$ on $S^{n-k-1} \times B^{k+1} =
\widetilde{U}(\rho/2)$. Hence $(\widetilde{U}(2s),g_\rho)$ is
isometric to an open subset of a manifold of the form 
$S^{n-k-1} \times N$ equipped with a product metric 
$g^{\text{round}} + g_N$, where $N$ is compact. By Proposition
\ref{appl.sbundles} the squared eigenvalues of the Dirac operator on
this product manifold are greater than or equal to
$(n-k-1)^2/4$. Writing the Rayleigh quotient of $\eta \psi_\rho$ we
obtain
\begin{equation} \label{rayquo}
\frac{(n-k-1)^2}{4} 
\leq 
\frac{\int_{\widetilde{U}(2s)}
|D^{g_\rho} (\eta \psi_{\rho}) |^2   \,dv^{g_\rho}}
{\int_{\widetilde{U}(2s)} 
| \eta \psi_{\rho} |^2  \,dv^{g_\rho}}.
\end{equation}
Since $D^{g_\rho} \psi_{\rho} = 0$ we have  
$D^{g_\rho} (\eta \psi_{\rho}) = 
\grad^{g_{\rho}} \eta \cdot \psi_{\rho}$ so
\begin{equation} \label{eq.dp}
|D^{g_\rho} (\eta\psi_\rho) |^2
= 
|\grad^{g_{\rho}} \eta \cdot \psi_{\rho}|^2
= 
|d\eta|_{g_{\rho}}^2 |\psi_{\rho}|_{g_{\rho}}^2.
\end{equation}
By definition $d \eta$ is supported in 
$\widetilde{U}(2s)\setminus \widetilde{U}(s)$. 
On $\widetilde{M} \setminus \widetilde{U}(2\rho)$
we have $g_{\rho} = F^2 g$. 
Moreover, by Relation (\ref{deta}) and since $F = 1/r$ on the 
support of $d \eta$, we have
$$ 
|d \eta|_{g_{\rho}}^2 = r^2 |d \eta|_g^2 \leq \frac{4r^2}{s^2}
$$
and hence 
$$ 
|D^{g_\rho} (\eta \psi_{\rho}) |^2 
\leq 
\frac{4r^2}{s^2} |\psi_{\rho}|^2,
$$
Since  $g_{\rho} = r^{-2} g$ on 
$\widetilde{U}(2s)\setminus \widetilde{U}(s)$ we have 
$dv^{g_\rho}=r^{-n} \,dv^g$. Using equation (\ref{eq.dp}) 
it follows that
\begin{eqnarray} \label{dpsi2}
\int_{\widetilde{U}(2s)} 
|D^{g_\rho} (\eta \psi_{\rho}) |^2 \,dv^{g_\rho}
&\leq& 
\frac{4}{s^2} 
\int_{\widetilde{U}(2s) \setminus \widetilde{U}(s)}
r^{2+(n-1)-n}
|r^{-\frac{n-1}{2}} \psi_\rho|^2 \,dv^g \nonumber\\
&\leq& 
\frac{8}{s} 
\int_{\widetilde{U}(2s)\setminus \widetilde{U}(s)} 
|F^{\frac{n-1}{2}}\psi_\rho|^2 \,dv^g,
\end{eqnarray}
where we also use that $r \leq 2s$ on the domain of integration.
Since $\eta \in [0,1]$ on 
$\widetilde{U}(2s)\setminus \widetilde{U}(s)$, since $\eta = 1$ on 
$\widetilde{U}(s)$ and since $g_{\rho} = r^{-2} g$ on 
$\widetilde{U}(s)\setminus \widetilde{U}(2 \rho)$, we have
\begin{eqnarray}
\int_{\widetilde{U}(2s)} 
|\eta \psi_{\rho} |^2 \,dv^{g_\rho} 
&\geq& 
\int_{\widetilde{U}(s)\setminus \widetilde{U}(2 \rho)} 
|\psi_{\rho} |^2 \,dv^{g_\rho} \nonumber\\
&=&
\int_{\widetilde{U}(s)\setminus \widetilde{U}(2 \rho)} 
r^{(n-1)-n}
| r^{-\frac{n-1}{2}}\psi_{\rho} |^2 \,dv^g \nonumber\\
&\geq&
\frac{1}{s}  
\int_{\widetilde{U}(s)\setminus \widetilde{U}(2 \rho)} 
|F^{\frac{n-1}{2}}\psi_{\rho}|_g^2 \,dv^g, \label{psi}
\end{eqnarray}
where we use that $r\leq s$ in the last inequality.
Plugging (\ref{dpsi2}) and (\ref{psi}) into (\ref{rayquo}) we get
$$
\frac{(n-k-1)^2}{4} 
\leq 
\frac{\frac{8}{s} 
\int_{\widetilde{U}(2s)\setminus \widetilde{U}(s)} 
|F^{\frac{n-1}{2}}\psi_\rho|^2 \,dv^g}
{\frac{1}{s}  
\int_{\widetilde{U}(s)\setminus \widetilde{U}(2 \rho)} 
|F^{\frac{n-1}{2}}\psi_{\rho}|_g^2 \,dv^g}
$$
and hence Proposition \ref{estim1} follows.
\end{proof}

\begin{proof}[Proof of Theorem \ref{main}]
As explained above we need to prove Relation (\ref{eqmain}), for a
contradiction assume that it is false. Then there is a strictly
decreasing sequence $\rho_i \to 0$ such that 
$\dim \ker D^g < \dim \ker D^{g_{\rho_i}}$ for all $i$. 
To simplify the notation for subsequences we define  
$E = \{ \rho_i : i\in\mN \}$. We have $0 \in \overline{E}$
and passing to a subsequence of $\rho_i$ means passing to a subset 
$E'\subset E$ of with $0 \in \overline{E'}$. 

Let  $m = \dim \ker D^g + 1$. For all $\rho \in E$ we can find 
$D^{g_{\rho}}$-harmonic spinors $\psi_{\rho}^1, \dots, \psi_{\rho}^m$
on $(\widetilde{M},g_\rho)$ such that 
\begin{equation} \label{normalis}
\int_{M\setminus U(s)}
\< \psi_{\rho}^j, \psi_{\rho}^k \> \,dv^g
=
\int_{\widetilde{M} \setminus \widetilde{U}(s)}
\< \psi_{\rho}^j,\psi_{\rho}^k \>  \,dv^g
=
\begin{cases}
1, &\text{if $j=k$;} \\
0, &\text{if $j\neq k$,}
\end{cases}
\end{equation}
where  $s \leq r_0< r_1/2$ is fixed as above. 
Let $\phi_{\rho}^j = F^{\frac{n-1}{2}}  \psi_{\rho}^j$. 
These spinor fields are defined on $M \setminus U(2 \rho)$ and by 
(\ref{confD}) they are $D^g$-harmonic.

\begin{step}
Let $\delta \in (0,R_{\text{max}})$. 
For $\rho>0$ small enough we have 
\begin{equation} \label{estim2}
\int_{M \setminus U(\delta)} | \phi_{\rho}^j|^2 \,dv^g 
\leq 
\frac{(n-k-1)^2 + 32}{(n-k-1)^2}.
\end{equation} 
\end{step}
By Proposition \ref{estim1} we have
$$
\int_{U(s)\setminus U(2\rho)} | \phi_{\rho}^j|^2 \,dv^g
\leq 
\frac{32}{(n-k-1)^2} \int_{U(2s)\setminus U(s)} 
|\phi_{\rho}^j|^2 \,dv^g. 
$$
and hence if $2\rho \leq \delta$ it follows that 
$$
\int_{U(s)\setminus U(\delta)} |\phi_{\rho}^j|^2 \,dv^g
\leq 
\frac{32}{(n-k-1)^2}   \int_{M \setminus U(s)} 
|\phi_{\rho}^j|^2 \,dv^g. 
$$
It follows that 
\begin{equation*}
\begin{split}
\int_{M \setminus U(\delta)} |\phi_{\rho}^j|^2 \,dv^g 
&=
\int_{M \setminus U(s)} | \phi_{\rho}^j|^2 \,dv^g
+ \int_{U(s) \setminus U(\delta)} 
|\phi_{\rho}^j|^2 \,dv^g \\
&\leq   
(1 + \frac{32}{(n-k-1)^2}) \int_{M \setminus U(s)} 
|\phi_{\rho}^j|^2 \,dv^g.
\end{split}
\end{equation*}
From (\ref{normalis}) we now obtain Inequality (\ref{estim2}).

\begin{step} 
There exists $E'\subset E$ with $0\in \overline{E'}$ and spinors  
$\Phi^1,\dots,\Phi^m \in C^1(M \setminus S)$, $D^g$-harmonic on 
$(M \setminus S,g)$ such that $\phi_{\rho}^j$ tend to $\Phi^j$ 
in $C^1_{\rm{loc}}(M\setminus S)$ as $\rho\to 0$, $\rho\in E'$.
\end{step}
Let $Z\in\mN$ be an integer, $Z>1/s$. By (\ref{estim2}) the sequence  
$\{ \phi_{\rho}^j \}_{\rho\in E}$ is bounded in 
$L^2(M \setminus U(1/Z))$. By Lemma \ref{lem.reg} it follows that 
$\{ \phi_{\rho}^j \}_{\rho\in E}$ is bounded in
$C^2(M \setminus U(2/Z))$ for all sufficiently large $Z$. 
For a fixed $Z_0>1/s$ we apply Lemma \ref{lem.ascoli}
and conclude that for any $j$ there is a subsequence 
$\{ \phi_{\rho}^j \}_{\rho\in E_0}$ of
$\{ \phi_{\rho}^j \}_{\rho\in E}$ that
converges in $C^1(M \setminus U(2/Z_0))$ to a spinor $\Phi_0^j$.
Similarly we construct further and further subsequences
$\{ \phi_{\rho}^j \}_{\rho\in E_i}$ converging to $\Phi_i^j$ 
in $C^1(M \setminus U(2/(Z_0+i)))$ with 
$E_i \subset E_{i-1} \subset \dots \subset E_0\subset E$, 
$0\in\overline{E_i}$. Obviously $\Phi_i^j$ extends $\Phi_{i-1}^j$. 
Define $E'\subset E$ as consisting of one $\rho_i$ from each $E_i$ 
chosen so that $\rho_i \to 0$ as $i \to \infty$. Then the sequence 
$\{ \phi_{\rho}^j \}_{\rho\in E'}$ converges in 
$C^1_{\text{loc}}(M \setminus S)$ to a spinor $\Phi^j$. 
As $\psi_{\rho}^j$ is $D^g$-harmonic on $(M\setminus U(2\rho))$ 
the $C^1_{\text{loc}}(M \setminus S)$-convergence implies that 
$D^g \Phi^j=0$ on $M\setminus S$. We have proved Step 2.

\begin{step}
Conclusion.
\end{step}
Let $j \in \{ 1, \dots,m\}$. 
By (\ref{estim2}) we conclude that 
$$
\int_{M \setminus S} |\Phi^j|^2 \,dv^g 
\leq 
\frac{(n-k-1)^2+32}{(n-k-1)^2}
$$
and hence $\Phi^j \in L^2(M)$. 
By Lemma \ref{lem.sing.rem} and elliptic regularity $\Phi^j$ is 
harmonic and smooth on all of $(M,g)$. 
Since $M \setminus U(s)$ is a relatively compact 
subset of $M \setminus S$ the normalization (\ref{normalis}) is 
preserved in the limit $\rho\to 0$ and hence 
$$
\int_{M \setminus U(s) } \<\Phi^j,\Phi^k\> \,dv^g = 
\begin{cases}
1, &\text{if $j=k$;} \\
0, &\text{if $j\neq k$.}
\end{cases}
$$
This proves that $\Phi^1, \dots, \Phi^m$ are linearly independent
harmonic spinors on $(M,g)$ and hence $\dim \ker D^g \geq m$ which
contradicts the definition of $m$. This proves Relation (\ref{eqmain})
and Theorem \ref{main}.
\end{proof}

\section{Proof of Theorem \ref{main_generic}}

The proof will follow the argument of \cite{baer_dahl02} so we
introduce notation in accordance to that paper. For a compact spin
manifold $M$ the space of smooth Riemannian metrics on $M$ is denoted
by ${\mathcal R}(M)$ and the subset of $D$-minimal metrics is
denoted by ${\mathcal R}_{\text{min}}(M)$.

From standard results in perturbation theory it follows that 
${\mathcal R}_{\text{min}}(M)$ is open in the $C^1$-topology on 
${\mathcal R}(M)$ and if ${\mathcal R}_{\text{min}}(M)$ is not empty then
it is dense in ${\mathcal R}(M)$ in all $C^k$-topologies, $k \geq 1$,
see for example \cite[Prop. 3.1]{maier97}. We define the word generic
to mean these open and dense properties satisfied by 
${\mathcal R}_{\text{min}}(M)$ if non-empty. Theorem \ref{main_generic}
is then equivalent to the following.
\begin{theorem} \label{reformulation_generic}
Let $M$ be a compact connected spin manifold. 
Then there is a $D$-minimal metric on $M$.
\end{theorem}
Before we start the proof we note the following consequence of Theorem
\ref{main}.
\begin{prop} \label{Dminaftersurgery}
Let $N$ be a compact spin manifold which has a $D$-minimal metric and
suppose that $M$ is obtained from $N$ by surgery of codimension $\geq 2$.
Then $M$ has a $D$-minimal metric. 
\end{prop}
\begin{proof}
This follows from Theorem \ref{main} since the left hand side of
(\ref{boundfromAS}) is the same for $M$ and $N$ while the right hand
side may only decrease. 
\end{proof}

From the proof of handle decompositions of bordisms we have the following.
\begin{prop} \label{bordism->surgery}
Suppose that $M$ is connected, $\dim M \geq 3$, and that $M$ is
spin bordant to a manifold $N$. Then $M$ can be obtained
from $N$ by a sequence of surgeries of codimension $\geq 2$. 
\end{prop}

\begin{proof}
The statement follows from \cite[VII Theorem 3]{kirby89} if $\dim M=3$.
If $\dim M\geq 4$, then we can do surgery in dimension $0$ and $1$ at a given 
spin cobordism between $M$ and $N$, and obtain
a connected, simply connected spin cobordism $W$ between $M$ and $N$. 
It then follows from \cite[VIII 3.1]{kosinski93} that one can obtain $M$ from
$N$ by surgeries of dimension $0,\ldots,n-2$.
\end{proof} 

\begin{proof}[Proof of Theorem \ref{reformulation_generic}]
From the solution of the Gromov-Lawson conjecture by Stolz
\cite{stolz92} together with knowledge of some explicit manifolds with
$D$-minimal metrics one can show that any compact spin manifold is
spin bordant to a manifold with a $D$-minimal metric, this is worked
out in detail in \cite[Prop. 3.9]{baer_dahl02}. We may thus assume
that the given manifold $M$ is spin bordant to a manifold $N$ equipped
with a $D$-minimal metric. The Theorem now follows from Propositions
\ref{Dminaftersurgery} and \ref{bordism->surgery} if $\dim M\geq 3$.

Now, let $\dim M=2$.
If $\alpha(M)=0$, then $M$ can 
be obtained by adding handles to $S^2$, i.e.\ by $0$-dimensional surgery.
If $\alpha(M)\neq 0$, then  $M$ can 
be obtained by adding handles to $T^2$ where $T^2$ carries the spin structure
with $\al\neq 0$. Any metric on $T^2$ with that spin structure has a 
$2$-dimensional kernel, and is thus $D$-minimal. 
With Proposition~\ref{Dminaftersurgery} we get 
Theorem~\ref{reformulation_generic} in the $2$-dimensional case.
\end{proof}

\providecommand{\bysame}{\leavevmode\hbox to3em{\hrulefill}\thinspace}
\providecommand{\href}[2]{#2}


\begin{thebibliography}{10}

\bibitem{adams75}
R.~A. Adams, \emph{Sobolev spaces}, Academic Press [A subsidiary of Harcourt
  Brace Jovanovich, Publishers], New York-London, 1975, Pure and Applied
  Mathematics, Vol. 65.

\bibitem{baer_dahl02}
C.~B{\"a}r and M.~Dahl, \emph{Surgery and the spectrum of the {D}irac
  operator}, J. Reine Angew. Math. \textbf{552} (2002), 53--76.

\bibitem{baer.schmutz:92}
C.~B{\"a}r and P.~Schmutz, \emph{{Harmonic spinors on Riemann surfaces.}}, Ann.
  Global Anal. Geom. \textbf{10} (1992), 263--273.

\bibitem{baum81}
H.~Baum, \emph{{S}pin-{S}trukturen und {D}irac-{O}peratoren {\"u}ber
  pseudoriemannschen {M}annigfaltigkeiten}, Teubner Verlag, 1981.

\bibitem{bourguignon_gauduchon92}
J.~P. Bourguignon and P.~Gauduchon, \emph{Spineurs, op{\'e}rateurs de {D}irac
  et variations de m{\'e}triques}, Math. Ann. \textbf{144} (1992), 581--599.

\bibitem{friedrich00}
T.~Friedrich, \emph{{D}irac {O}perators in {R}iemannian {G}eometry}, {G}raduate
  {S}tudies in {M}athematics 25, AMS, Providence, Rhode Island, 2000.

\bibitem{gilbarg.trudinger77}
D.~Gilbarg and N.~Trudinger, \emph{Elliptic partial differential equations of
  second order}, Grundlehren der mathematischen Wissenschaften, no. 224,
  Springer Verlag, 1977.

\bibitem{hijazi01}
O.~Hijazi, \emph{{Spectral properties of the {D}irac operator and geometrical
  structures.}}, {Ocampo, Hernan (ed.) et al., Geometric methods for quantum
  field theory. Proceedings of the summer school, Villa de Leyva, Colombia,
  July 12-30, 1999. Singapore: World Scientific. 116-169 }, 2001.

\bibitem{hitchin74}
N.~Hitchin, \emph{Harmonic spinors}, Advances in Math. \textbf{14} (1974),
  1--55.

\bibitem{kirby89}
R.~C. Kirby, \emph{The topology of {$4$}-manifolds}, Lecture Notes in
  Mathematics, vol. 1374, Springer-Verlag, Berlin, 1989.

\bibitem{kosinski93}
A.~A. Kosinski, \emph{Differential manifolds}, Pure and Applied Mathematics,
  vol. 138, Academic Press Inc., Boston, MA, 1993.


\bibitem{lawson.michelsohn89}
H.~B. Lawson and M.-L. Michelsohn, \emph{Spin geometry}, Princeton Mathematical
  Series, vol.~38, Princeton University Press, Princeton, NJ, 1989.

\bibitem{maier97}
S.~Maier, \emph{Generic metrics and connections on {S}pin- and {S}pin{$\sp
  c$}-manifolds}, Comm. Math. Phys. \textbf{188} (1997), no.~2, 407--437.

\bibitem{stolz92}
S.~Stolz, \emph{Simply connected manifolds of positive scalar curvature}, Ann.
  of Math. (2) \textbf{136} (1992), no.~3, 511--540.


\end{thebibliography}

\end{document}